\documentclass[12pt,a4paper]{article}
\usepackage{amssymb}
\usepackage{graphicx,amssymb,amsfonts,epsfig,amsthm,a4,amsmath,url}
\usepackage[latin1]{inputenc}

\newtheorem{thm}{Theorem}[section]
\newtheorem{prop}[thm]{Proposition}
\newtheorem{lem}[thm]{Lemma}
\newtheorem{cor}[thm]{Corollary}
\theoremstyle{definition}
\newtheorem{defn}[thm]{Definition}

\newtheorem*{Ack}{Acknowledgements}
\newtheorem{rem}[thm]{Remark}
\newtheorem{exe}[thm]{Example}
\newtheorem{conj}{Conjecture}
\newtheorem{que}{Question}

\newcommand{\BP}{\text{BP}}
\newcommand{\PHO}{$\text{(BP}_0\text{)}$}
\newcommand{\FV}{\text{FH}}
\newcommand{\FHO}{$\text{(FH}_0\text{)}$}
\newcommand{\Z}{\mathbf{Z}}

\newcommand{\N}{\mathbf{N}}
\newcommand{\R}{\mathbf{R}}
\newcommand{\C}{\mathbf{C}}

\newcommand{\D}{\mathbf{D}}
\newcommand{\Ker}{\text{Ker}}

\newcommand{\HH}{\mathcal{H}}

\newcommand{\eps}{\varepsilon}

\newcommand{\LL}{\mathcal{L}}

\title{Isometric group actions on Banach spaces and representations vanishing at infinity}
\author{Yves de Cornulier, Romain Tessera, Alain Valette}
\date{\today}

\begin{document}

\baselineskip=16pt

\maketitle

\begin{abstract}
Our main result is that the simple Lie group $G=Sp(n,1)$ acts
properly isometrically on $L^p(G)$ if $p>4n+2$. To prove this, we
introduce property $({\BP}_0^V)$, for $V$ be a Banach space: a
 locally compact group $G$ has property $({\BP}_0^V)$ if every affine
isometric action of $G$ on $V$, such that the linear part is a
$C_0$-representation of $G$, either has a fixed point or is
metrically proper. We prove that solvable groups, connected Lie
groups, and linear algebraic groups over a local field of
characteristic zero, have property $({\BP}_0^V)$. As a consequence
for unitary representations, we characterize those groups in the
latter classes for which the first cohomology with respect to the
left regular representation on $L^2(G)$ is non-zero; and we
characterize uniform lattices in those groups for which the first
$L^2$-Betti number is non-zero.
\medskip

    \hfill\break
\noindent {\sl Mathematics Subject Classification:} Primary 22D10;
Secondary 20J06, 22D10, 22E40, 43A15. \hfill\break {\sl Key words
and Phrases:} Affine isometries, isometric representations,
1-cohomology, vanishing of coefficients.
\end{abstract}


\section{Introduction}

The study of affine isometric actions on Hilbert spaces has, in
recent years, found applications ranging from geometric group
theory, to rigidity theory, to $K$-theory of $C^*$-algebras. For
instance, if $G$ is $\sigma$-compact, $G$ has Kazhdan's Property (T)
if and only if every affine isometric action on a Hilbert space has
a fixed point (see \cite[Chapter~2]{BHV}). This is known to have
strong group-theoretic consequences on~$G$: for instance, this
implies that $G$ is compactly generated and has compact
abelianization (see \cite[Chap. 2]{BHV} for a direct proof). On the
other hand, a group is said to be a-T-menable if it admits a proper
affine action on a Hilbert space, and a deep result of Higson and
Kasparov \cite{HiKa} asserts that an a-T-menable group satisfies the
Baum-Connes conjecture.

These notions naturally extend to affine isometric actions on Banach
spaces. Generalizations of Property (T) to the context of uniformly
convex Banach spaces have proven to be fruitful, producing new
rigidity results \cite{BFGM,FM}; on the other hand, a result of
Kasparov and G. Yu \cite{KaYu} asserts that the Novikov conjecture
holds for finitely generated groups embedding uniformly in
superreflexive Banach spaces (in particular for those admitting a
proper isometric action on such a space).

Guoliang Yu recently proved \cite{Yu} that every Gromov hyperbolic group
admits a proper isometric action on the uniformly convex space
$\ell^p( \Gamma\times\Gamma)$, for $p\gg 0$. This result contrasts
with the existence of Gromov hyperbolic groups with Property (T).

Our main result is the following theorem.

\begin{thm}\label{thm1}
Let $k$ be a local field. Let $G$ be a simple algebraic group of
rank $1$ over $k$. Let $p_0$ be the Hausdorff dimension of the
visual boundary of $G$. Then, for every $p>\max\{1,p_0\}$, there
exists a proper affine action of $G$ on $L^p(G)$ with linear part
$\lambda_{G,p}$.
\end{thm}

Note that this result cannot be extended to simple groups of rank
$\geq 2$ (see \cite[Theorem~B]{BFGM}). We also prove (see
Proposition \ref{walls}) that every countable a-T-menable group
admits a proper affine action on some $L^p$-space, for every $p\geq
1$.

\medskip

A key ingredient of the proof of Theorem \ref{thm1} is a result of
Pansu~\cite{Pan,Pansu}, who computes the first $L^p$-cohomology for
semi-simple Lie groups for $1< p<\infty$. The first $L^p$-cohomology
actually coincides (see Section~\ref{sectioncohLp}) with the first
cohomology of the group with values in the right regular
representation on $L^p(G)$. The proof of Theorem~\ref{thm1} then
consists in showing that non-trivial $1$-cocycles on such
representations are automatically proper. This latter fact is part
of a more general phenomenon: the properness of non-trivial
$1$-cocycles on an isometric $L^p$-representation $\pi$ of a group
$G$ is actually true under very general assumptions on $G$ and
$\pi$.

Our approach was initially motivated by the following example. The
cyclic group $\Z$ acts naturally on $\ell^2(\Z)$; the corresponding
operator $T$, given by the action of the positive generator of $\Z$
is usually called the bilateral shift. Now take $f\in\ell^2(\Z)$,
and let us consider the affine isometry $T_f$ of $\ell^2(\Z)$ given
by $T_f(v)=Tv+f$. It is immediately checked that this isometry has a
fixed point if and only if $f\in\textnormal{Im}(T-1)$. We show that
otherwise the corresponding action is proper, that is, for
every/some $v\in\ell^2(\Z)$, $\|{T_f}^n(v)\|\to\infty$ when
$|n|\to\infty$. Our aim is to make this observation systematic.

One essential feature in the above context is that the
representation of $\Z$ on $\ell^2(\Z)$ is $C_0$. In a general
context, let $V$ be a Banach space. An isometric linear
representation $\pi$ of a locally compact group $G$ is $C_0$ if for
every $L\in V^*$ (the topological dual), and every $v\in V$, we have
$L(\pi(g)v)\to 0$ when $g$ tends to infinity. In other words,
$\pi(g)v$ weakly tends to zero for every $v\in V$.

\begin{defn} Let $V$ be a Banach space.
A locally compact group $G$ has {\it Property $(\BP_0^V)$} if, for
every $C_{0}$ isometric linear representation $\pi$ of $G$ on $V$,
any affine isometric action of $G$ with linear part $\pi$ either has
a bounded orbit or is proper. We say that $G$ has property {\PHO} if
it has $(\BP_0^V)$ for every Hilbert space $V$.
\end{defn}

The acronym $({\BP}_0)$ stands for ``Bounded or Proper with respect
to $C_0$-repre\-sentations".

Thus the observation above amounts to prove that $\Z$ has Property
$({\BP}_0)$. This is part of the following result.

\begin{prop} (see Proposition \ref{PH0heredit})
Let $G$ be a locally compact group, and $V$ a Banach space.

(1) Suppose that $G$ has two non-compact normal subgroups
centralizing each other. Then $G$ has Property $({\BP}_0^V)$.

(2) Suppose that $G$ has a non-compact normal subgroup with Property
$({\BP}_0^V)$. Then $G$ has Property $({\BP}_0^V)$.
\end{prop}

\begin{cor} (see Corollary \ref{solvable} and Proposition \ref{groupswithPH0})
Let $V$ be a Banach space.
\begin{enumerate}
\item [1)] Every solvable locally compact group has Property $({\BP}_0^V)$.
\item [2)] Every connected Lie group or linear algebraic group over a $p$-adic field
has property $({\BP}_0^V)$.
\end{enumerate}\label{CorIntro}
\end{cor}

As an application of Corollary \ref{CorIntro}, we characterize in
Proposition \ref{H1reg} those connected Lie groups and linear
algebraic groups over a $p$-adic field, such that the first
cohomology of $G$ with coefficients in the left regular
representation $\lambda_G$ on $L^2(G)$ is non-zero; this generalizes
a result of Guichardet (Proposition 8.5 in Chapter III of
\cite{Guich}) for simple Lie groups.

\begin{prop} (see Proposition \ref{H1reg}) Let $G$ be a connected Lie
group or $G=\mathbf{G}(\mathbf{K})$,
the group of $K$-points of a linear algebraic group $\mathbf{G}$
over a local field $\mathbf{K}$ of characteristic zero. Assume $G$
non-compact. Then the following are equivalent\begin{itemize}
    \item[(i)] $H^1(G,\lambda_G)\neq 0$.
    \item[(ii)] Either $G$ is amenable, or there exists a compact normal subgroup $K\subset G$ such
    that $G/K$ is isomorphic to $\textnormal{PSL}_2(\R)$ (case of Lie groups),
    or a simple algebraic group of rank one (case of an algebraic group over a $p$-adic field).
    \end{itemize}
\end{prop}

We also characterize those uniform lattices $\Gamma$ in a group as
above, whose first $L^2$-Betti number $\beta^1_{(2)}(\Gamma)$ is
non-zero. For uniform lattices in connected Lie groups, this gives a
new proof of Theorem 4.1 in \cite{Eck}.

\begin{cor} (see Corollary \ref{bettiL2}) Let $G$ be a connected
Lie group or $G=\mathbf{G}(\mathbf{K})$
where $\mathbf{K}$ is a local field of characteristic zero; let
$\Gamma$ be a uniform lattice in $G$. If the first $L^2$-Betti
number $\beta^1_{(2)}(\Gamma)$ is non-zero, then $\Gamma$ is
commensurable either to a non-abelian free group or to a surface
group.
\end{cor}

In contrast with property $({\BP}_0)$, we have

\begin{prop}(see Proposition \ref{prop:affine_Z_non_wc0})
There exists an affine isometric action of $\Z$ on a complex Hilbert
space, that is neither proper nor bounded, and whose linear part has
no finite-dimensional subrepresentation.
\end{prop}

This result can be extended to $\R$ in view of the following result.

\begin{prop}(see Proposition \ref{Isomexpon})
Every isometric action of $\Z$ on a complex Hilbert space can be
extended to a continuous action of~$\R$.
\end{prop}

While property $({\BP}_0)$ is a rule for connected Lie groups or
linear algebraic groups over local fields of characteristic zero,
this is certainly not the case for discrete groups:

\begin{prop} (see Proposition \ref{amalgams} and Corollary \ref{surfgrps})
Non-abelian free groups and surface groups do not have property
$({\BP}_0)$.\label{nonBP0}
\end{prop}

The proof of Proposition \ref{nonBP0} follows from a direct and
simple construction for free groups. However, for surface groups,
our proof relies on von Neumann algebra arguments as well as the
analytical zero divisor Conjecture.

\begin{Ack}
We thank Emmanuel Breuillard for his contribution to the proof of
Proposition \ref{prop:affine_Z_non_wc0}, and Uffe Haagerup for
pointing out reference \cite{Aag}.
\end{Ack}

\section{Actions on Banach spaces}\label{PHOsection}

We define a {\it Banach pair} as a pair $(V,\LL)$ where $V$ is a
Banach space and $\LL$ is a linear subspace of $V^*$. We call it a
Banach-Steinhaus pair if it satisfies the Banach-Steinhaus Property:
any subset $X\subset V$ is bounded if and only if $L(X)$ is bounded
for every $L\in\LL$. For instance, the Banach-Steinhaus Theorem
states that $(V^*,V)$ is a Banach-Steinhaus pair for every Banach
space $V$, and in particular $(V,V^*)$ is a Banach-Steinhaus pair.
If $(V,\LL)$ is a Banach-Steinhaus pair, and if $W$ is a closed
subspace of $V$, then $(W,\LL|_W)$ is clearly a Banach-Steinhaus
pair, where $\LL|_W$ is the set of restrictions of $L\in\LL$ to $W$.

If $(V,\LL)$ is a Banach pair, we say that an isometric linear
action $\pi$ of a locally compact group $G$ on $V$ is $C_0^\LL$ if
$L(\pi(g)v)$ tends to zero when $g\to\infty$, for every $v\in V$ and
$L\in\LL$. Note that being $C_0^\LL$ definitely depends on $\LL$
(see the example below); however when the context is clear we write
it $C_0$.

\begin{exe}
Let $G$ be a discrete, infinite group. Consider its regular
$\ell^1$-representation. Then it is $C_0^{c_0(G)}$ but not
$C_0^{\ell^{\infty}(G)}$. Note that both are Banach-Steinhaus pairs.
This example motivates the introduction of Banach-Steinhaus pairs
different from $(V,V^*)$.
\end{exe}

If $\pi$ is a $C_0$ representation as above, and if $W$ is a closed
invariant subspace, defining a subrepresentation $\pi|^W$, then
$\pi|^W$ is $C_0^{\LL|_W}$.

\begin{defn} Let $(V,\LL)$ be a Banach pair.
A locally compact group $G$ has {\it Property $(\FV_0^{(V,\LL)})$
(respectively $(\BP_0^{(V,\LL)})$)} if, for every
$C_{0}$-representation $\pi$ of $G$, any affine isometric action of
$G$ on $V$ with linear part $\pi$ has a bounded orbit (resp. either
has a bounded orbit or is proper).

We say that $G$ has {\it Property $(\FV_0^{([V],\LL)})$} if it has
Property $(\FV_0^{(W,\LL|_W)})$ for every closed subspace $W$ of
$V$. We define analogously Property $(\BP_0^{([V],\LL)})$.

Similarly, we say that $G$ has {\it Property $(\FV_0^{V})$
(respectively $(\BP_0^{V})$)} if it has {\it Property
$(\FV_0^{(V,\LL)})$ (respectively $(\BP_0^{(V,\LL)})$)} for
$\LL=V^*$.
\end{defn}

When the space $V$ is superreflexive, i.e. isomorphic to a uniformly
convex space, it is known that every nonempty bounded subset has a
unique circumcenter (also called Chebyshev center, see p. 27 in
\cite{BL}). As a consequence, every isometric action with a bounded
orbit on $V$ has a globally fixed point.

\begin{lem}
Let a compact group $K$ act by affine isometries on a Banach space.
Then it fixes a point.\label{lem:compact_fixes}
\end{lem}
\begin{proof}Let $\Omega$ be an orbit. As $\Omega$
is compact, its closed convex hull $X$ is also compact (see for
instance \cite[Theorem~3.25]{Rudin}). As $K$ is amenable and acts on
$X$ by affine transformations, it has a fixed point.\end{proof}

\begin{lem}\label{compactnormal} Let $(V,\LL)$ be a Banach pair.
    Let $K$ be a compact, normal subgroup of $G$. The following are
equivalent.
    \begin{enumerate}
        \item  [(i)] $G$ has Property $(\FV_0^{([V],\LL)})$
(resp. $(\BP_0^{([V],\LL)})$);

        \item  [(ii)] $G/K$ has Property $(\FV_0^{([V],\LL)})$
(resp. $(\BP_0^{([V],\LL)})$).
    \end{enumerate}
\end{lem}
\begin{proof}The implication (i)$\Rightarrow$(ii) is clear. Suppose
that $G/K$ has Property $(\BP_0^{([V],\LL)})$. By Lemma
\ref{lem:compact_fixes}, the set $W$ of $K$-fixed points is a
non-empty closed affine subspace; moreover it is $G$-invariant. As
$G$ has Property $(\BP_0^{(W,\LL|_W)})$, its action on $W$, and
therefore on $V$, is either bounded or proper. The case of Property
$(\FV_0^{([V],\LL)})$ is similar.\end{proof}

\begin{lem}\label{C0+subgroups}Let $(V,\LL)$ be a Banach-Steinhaus pair.
Let $H,K$ be closed, non-compact subgroups of the locally compact
group $G$ which centralize each other. Let $\alpha$ be an affine
isometric action of $G$ on $V$, whose linear part $\pi$ is a
$C_{0}$-representation. Then either $\alpha|_H$ and $\alpha|_K$ are
both bounded, or they are both proper.\end{lem}

\begin{proof} Set $b(g)=\alpha(g)(0)$. We assume that $\alpha|_H$ is not
proper, i.e.
$$M=:\liminf_{h\rightarrow\infty,h\in H}\|b(h)\|\,<\,\infty.$$
For $k\in K,\,h\in H$, the 1-cocycle relation gives
$$\pi(k)b(h)+b(k)\,=\,b(kh)\,=\,b(hk)\,=\,\pi(h)b(k)+b(h),$$
which we will use in the following form:
$$b(k)\,=\pi(h)b(k) + (1-\pi(k))b(h).$$
Then, for every $L\in\LL$ we have
$$L(b(k))=L(\pi(h)b(k)) + L((1-\pi(k))b(h)),$$ and thus
$$|L(b(k))|\le |L(\pi(h)b(k))| + |L((1-\pi(k))b(h))|\le
|L(\pi(h)b(k))| + 2\|L\|\|b(h)\|.$$ Taking the inferior limit when
$h\to\infty$ in $H$, we obtain
$$|L(b(k))|\le 2\|L\|M.$$
Thus $L(b(K))$ is bounded for every $L$; as $(V,\LL)$ is a
Banach-Steinhaus pair this means that $b(K)$ is bounded.

Inverting the roles of $H$ and $K$, we can easily conclude.
\end{proof}

The following proposition is an immediate consequence of Lemma
\ref{C0+subgroups}, by taking $H=G$ and $K=Z(G)$.

\begin{prop}\label{noncompactcentre}
    Let $G$ be a locally compact group with non-compact centre (e.g.
    a non-compact, locally compact abelian group). Then $G$ has property
    $(\BP_0^{(V,\LL)})$ for every Banach-Steinhaus pair $(V,\LL)$.\qed
\end{prop}

In order to enlarge the class of groups for which we are able to
prove Property $({\BP}_0)$, we need the following lemma.

\begin{lem}\label{normalsubgr}
    Let $(V,\LL)$ be a Banach-Steinhaus pair.
Let $\alpha$ be an affine isometric action of $G$ on $V$, with
linear part a $C_{0}$-representation $\pi$. Set $b(g)=\alpha(g)(0)$.
Let $H$ be a closed, non-compact subgroup of $G$. Assume that there
exists a sequence $(g_{k})_{k\geq 1}$ in $G$, going to infinity,
such that
\begin{itemize}
    \item  the sequence $(b(g_k))$ is bounded in $V$;
    \item  for every $h\in H$, the set $\{g_{k}^{-1}hg_{k}|\,k\ge 1\}$
is relatively compact in $G$.
\end{itemize}
Then $\alpha|_{H}$ is bounded.
\end{lem}

\begin{proof} Fix $M>0$ such that $\|b(g_{k})\|\leq M$
    for every $k\geq 1$, and, for $h\in H$ define $K_h$ as the closure of
the set $\{g_{k}^{-1}hg_{k}|\,k\ge 1\}$, which is compact by
assumption. Let us show that $\|b(h)\|\leq 2M$ for every $h\in H$.
Noting that $hg_k=g_kh_k$, where $h_k=g_{k}^{-1}hg_{k}$, we expand
$b(hg_{k})-b(g_{k})-b(h)$ in two ways, and we obtain
$$\pi(g_{k})b(h_{k})-b(h)=(\pi(h)-1)b(g_{k}).$$
So, for every $L\in\LL$,
$$|L(b(h))|\le |L(\pi(g_{k})b(h_{k}))|+|L((\pi(h)-1)b(g_{k}))|$$ $$\le
|L(\pi(g_{k})b(h_{k}))|+2\|L\|\|b(g_{k})\|.$$

Using the assumption that $h_{k}\in K_{h}$ for every $k$, and the
fact that for a $C_{0}$-representation the decay of coefficients to
$0$ is uniform on compact subsets of the ambient Banach space, we
get for $k\rightarrow\infty$,
$$|L(b(h))|\le 2M\|L\|.$$
As $(V,\LL)$ is a Banach-Steinhaus pair, this implies that $b(H)$ is
bounded.\end{proof}

The following lemma is a kind of a geometric Hahn-Banach statement.

\begin{lem}
Let $(V,\LL)$ be a Banach-Steinhaus pair. Then there exists $\eta>0$
with the following property: for every $v\in V$, there exists
$L\in\LL$ such that $\|L\|\le 1$ and
$L(v)\ge\eta\|v\|$.\label{lem:sep_conv}
\end{lem}
\begin{proof}
Suppose the contrary. For every $n$, there exists $v_n\in V$ of norm
one such that for every $L\in\LL$, we have $L(v_n)<2^{-n}\|L\|$. Set
$X=\{2^nv_n|n\ge 0\}$. Then $L(X)$ is bounded for every $L\in\LL$;
by the Banach-Steinhaus Property, $X$ is bounded; this is a
contradiction.
\end{proof}

\begin{lem}
Let $(V,\LL)$ be a Banach-Steinhaus pair. Let $G$ be a locally
compact group, and $N$ a non-compact, closed normal subgroup. Let
$\alpha$ be an affine isometric action of $G$ on $V$ whose linear
part is $C_0$.

\begin{itemize}
\item[(1)] Suppose that $\alpha|_N$ is
bounded. Then $\alpha$ is also bounded.

\item[(2)] Suppose that $\alpha|_N$ is
proper. Then $\alpha$ is also proper.
\end{itemize}\label{lem:bounded_inherited}
\end{lem}
\begin{proof}
(1) For $M\ge 0$, define $A_M$ as the set of all $x\in V$ whose
$N$-orbit has diameter at most $M$. Clearly $A_M$ is $G$-invariant.
By assumption, for some $M$ (which we fix now), the set $A_M$ is
non-empty. We claim that it is bounded, allowing us to conclude.

Consider $x,y\in A_M$, and set $v=x-y$. Then, for $h\in N$
$$\pi(h)v-v=\alpha(h)x-x-\alpha(h)y+y.$$
So
$$\|\pi(h)v-v\|\le\|\alpha(h)x-x\|+\|y-\alpha(h)y\|\le 2M.$$
Fix $\eta$ and $L$ as in Lemma \ref{lem:sep_conv}. Then
$$\eta\|v\|\le L(v)\le |L(v)-L(\pi(h)v)|+|L(\pi(h)v)|\le
2M+|L(\pi(h)v)|.$$ As $N$ is non-compact, letting $h\to\infty$, we
obtain $\|v\|\le 2M/\eta$. Thus the diameter of $A_M$ is bounded by
$2M/\eta$.

(2) Suppose by contradiction that $\alpha|_N$ is proper and $\alpha$
is not proper. Then there exists a sequence $(g_k)$ in $G$, tending
to infinity, such that $(b(g_k))$ is bounded. As $\alpha|_N$ is
unbounded, there exists, by Lemma \ref{normalsubgr}, an element
$n\in N$ such that the sequence $(g_{k}^{-1}ng_{k})_{k\geq 1}$ is
not relatively compact in $N$; extracting if necessary we can
suppose that it tends to infinity. Therefore, as $\alpha|_N$ is
proper, $\|b(g_{k}^{-1}ng_{k})\|$ tends to infinity. But it is
bounded by $2\|b(g_k)\|+\|b(n)\|$, which is bounded, a
contradiction.
\end{proof}

From this we deduce the following.

\begin{prop}\label{PH0heredit}
Let $(V,\LL)$ be a Banach-Steinhaus pair, and let $G$ be a locally
compact group. Let $N$ be a non-compact, closed, normal subgroup of
$G$. If either the centralizer $C_{G}(N)$ of $N$ is non-compact, or
$N$ has Property $(\BP_0^{(V,\LL)})$, then $G$ also has Property
$(\BP_0^{(V,\LL)})$.
\end{prop}

\begin{proof} Let $\alpha$ be an affine isometric
action of $G$, with linear part a $C_{0}$-representation $\pi$. If
$\alpha|_{N}$ is bounded, then $\alpha$ is bounded by Lemma
\ref{lem:bounded_inherited}(1). Assume then that $\alpha|_{N}$ is
unbounded. Then $\alpha|_{N}$ is proper (in case $C_{G}(N)$ is
non-compact, this follows from lemma \ref{C0+subgroups}).
Accordingly, by Lemma \ref{lem:bounded_inherited}(2), $\alpha$ is
proper.
\end{proof}

\begin{cor}
Let $(V,\LL)$ be a Banach-Steinhaus pair. Then Properties
$(\BP_0^{([V],\LL)})$ and $(\FV_0^{([V],\LL)})$ are preserved by
extensions.\label{cor:PHO_extensions}
\end{cor}
\begin{proof} Let $1\to N\to G\to H\to 1$ be an extension of locally
compact groups, and suppose that $N$ and $H$ have Property
$(\BP_0^{([V],\LL)})$. If $N$ is compact, then, by Lemma
\ref{compactnormal}, since $H$ has Property $(\BP_0^{([V],\LL)})$,
so does $G$. If $N$ is not compact, then, since it has Property
$(\BP_0^{([V],\LL)})$, by Proposition \ref{PH0heredit}, $G$ has
Property $(\BP_0^{([V],\LL)})$. The case of Property
$(\FV_0^{([V],\LL)})$ is similar (and easier).\end{proof}

\begin{cor}\label{solvable}
    Locally compact, solvable groups have Property
$(\BP_0^{(V,\LL)})$ for every Banach-Steinhaus pair $(V,\LL)$.
\label{cor:solvable_PH0}
\end{cor}
\begin{proof} Since, using Proposition \ref{noncompactcentre}, locally
compact abelian groups have Property $(\BP_0^{(V,\LL)})$, this
follows from Corollary \ref{cor:PHO_extensions}.\end{proof}

\begin{lem}
Let $(V,\LL)$ be a Banach pair. Property $(\BP_0^{(V,\LL)})$ is
inherited from cocompact subgroups.\label{lem:cocompact}
\end{lem}
\begin{proof} The proof is straightforward.\end{proof}

\begin{prop}\label{groupswithPH0}
    Connected Lie groups, and linear algebraic groups over $p$-adic
fields,
    have Property $(\BP_0^{(V,\LL)})$ for every Banach-Steinhaus pair $(V,\LL)$.
\end{prop}
\begin{proof} This follows from Lemma \ref{lem:cocompact} and Corollary
\ref{cor:solvable_PH0}, since $G$ contains a solvable cocompact
subgroup $H$: for linear algebraic groups over local fields of
characteristic zero, this follows from
\cite[Théorème~8.2]{BorelTits}; for Lie groups, taking the quotient
by the maximal solvable normal subgroup, we can also use
\cite[Théorème~8.2]{BorelTits}.\end{proof}

\section{Proper affine actions of rank 1 groups on $L^p$-spaces}

\subsection{Spaces with measured walls, and the non-archimedean case}

Recall that a locally compact $\sigma$-compact group is a-T-menable
if it acts properly isometrically on some Hilbert space.

\begin{prop}\label{walls} Let $\Gamma$ be a countable, discrete group.
The following are equivalent:
\begin{enumerate}
\item [i)] $\Gamma$ is a-T-menable;
\item [ii)] for every $p\geq 1$, the group $\Gamma$ acts properly isometrically on some $L^p$-space.
\end{enumerate}
\end{prop}

\begin{proof} We prove the non-trivial implication $(i)\Rightarrow (ii)$. We recall from
\cite{CMV} that a space with measured walls is a
4-tuple $(X,{\cal W,B},\mu)$ where $X$ is a set, ${\cal W}$ is a set
of partitions of $X$ into 2 classes (called walls), ${\cal B}$ is a
$\sigma$-algebra of sets on ${\cal W}$, and $\mu$ is a measure on
${\cal B}$ such that, for every pair $x,y$ of distinct points in
$X$, the set $\omega(x,y)$ of walls separating $x$ from $y$ belongs
to $\cal{B}$ and satisfies $w(x,y)=:\mu(\omega(x,y))<\infty$.

It was proved in Proposition 1 of \cite{CMV} that a countable group
is a-T-menable if and only if it admits a proper action on some
space with measured walls (by this we mean that $\Gamma$ preserves
the measured wall space structure, and that the function $g\mapsto
w(gx,x)$ is proper on $\Gamma$).

A half-space in a space with measured walls $X$ is a class of the
partition defined by some wall in ${\cal W}$. Let $\Omega$ be the
set of half-spaces, $p:\Omega\rightarrow {\cal W}$ the canonical map
(associating to any half-space the corresponding wall), ${\cal A}
=:p^{-1}({\cal B})$ the pulled-back $\sigma$-algebra, and $\nu$ the
pulled-back measure defined by
$$\nu(A)=\frac{1}{2}\int_{{\cal W}} card(A\cap p^{-1}(x))\,d\mu(x)$$
for $A\in{\cal A}$. Let $\chi_x$ be the characteristic function of
the set of half-spaces through $x$. For $x,y\in X$, we define a
function $c(x,y)\in L^p(\Omega,\nu)$ by:
$$c(x,y)=\chi_x-\chi_y.$$
Suppose that $\Gamma$ acts properly on $(X,{\cal W,B},\mu)$. For
$p\geq 1$, let $\pi_p$ denote the quasi-regular representation of
$\Gamma$ on $L^p(\Omega,\nu)$. Observe that:
\begin{itemize}
\item $c(x,y)+c(y,z)= c(x,z)$;
\item $c(gx,gy)=\pi_p(g)c(x,y)$;
\item $\|c(x,y)\|_p^p = w(x,y)$
\end{itemize}
for every $x,y,z\in X,\,g\in\Gamma$. Fixing a base-point $x_0\in X$,
the map
$$b:\Gamma\rightarrow L^p(\Omega,\nu):g\mapsto c(gx_0,x_0)$$
defines a 1-cocycle in $Z^1(\Gamma,\pi_p)$. Since $\|b(g)\|_p
=w(gx_0,x_0)^{1/p}$, this cocycle is proper, so that the
corresponding affine isometric action is proper.
\end{proof}

\begin{rem}: What the above proof really shows is that every locally
compact group acting properly on a space with
measured walls, admits a proper action on some $L^p$-space, for
every $p\geq 1$. Several non-discrete examples appear in \cite{CMV}.
\end{rem}

A tree $X=(V,E)$ is an example of a space with measured walls (with
${\cal W}=E,\,\mu=$ counting measure). The set $\Omega$ of
half-spaces identifies with the set $\mathbb{E}$ of oriented edges.
Suppose that a locally compact group $G$ acts properly cocompactly
on a tree $X$. We choose a reference edge $e_0\in \mathbb{E}$ and
use it to lift the cocycle $b\in Z^1(G,\ell^p(\mathbb{E}))$ from the
previous proof to a cocycle $\tilde{b}\in Z^1(G,L^p(G))$, by the
formula $(\tilde{b}(g))(h)=(b(g))(he_0)$. Then
$$\|\tilde{b}(g)\|_p^p
= \frac{m_0d(gx_0,x_0)}{k},$$ where $m_0$ is the Haar measure of the
stabilizer of $e_0$ in $G$, and $k$ is the number of orbits of $G$
in $\mathbb{E}$. This shows that $\tilde{b}$ is a proper cocycle. We
have proved:

\begin{prop}\label{trees} Let $G$ be a locally compact group. If $G$ acts
properly cocompactly on a tree (e.g. if $G$ is a rank 1 simple
algebraic group over a non-Archimedean local field), then for every
$p\geq 1$, the group $G$ admits a proper isometric action on
$L^p(G)$, with linear part the left regular representation
$\lambda_{G,p}$. \hfill$\qed$
\end{prop}

\subsection{The Lie group case}\label{sectioncohLp}

Let $M$ be a Riemannian manifold equipped with its Riemannian
measure $\mu$.  Fix $p>1.$ Denote by $D_p(M)$ the vector space of
differentiable functions whose gradient is in $L^p(TM)$. Equip
$D_p(M)$ with a pseudo-norm $\|f\|_{D_p}=\|\nabla f\|_p$, which
induces a norm on $D_p(M)$ modulo the constants. Denote by $\D_p(M)$
the completion of this normed vector space. We have
$W^{1,p}(M)=L^p(M)\cap D_p(M)$. Hence, $W^{1,p}(M)$ canonically
embeds in $\D_p(M)$ as a subspace.

The first $L^p$-cohomology of $M$ is the quotient space
$$H_p^1(M)=\D_p(M)/W^{1,p}(M).$$ The first {\it reduced}
$L^p$-cohomology of $M$ is the quotient space
$$\overline{H_p}^1(M)=\D_p(M)/\overline{W^{1,p}(M)},$$
where $\overline{W^{1,p}(M)}$ is the closure of $W^{1,p}(M)$ in the
Banach space $\D_p(M)$. Note that the two latter spaces coincide if
and only if the norm $\|\cdot\|_{D_p}$ on the Sobolev space
$W^{1,p}(M)$ is equivalent to the usual Sobolev norm
$\|\cdot\|_p+\|\cdot\|_{D_p}$, that is, if $M$ satisfies the strong
Sobolev inequality in $L^p$: $\|f\|_p\leq C\|\nabla f\|_{p}.$ If the
group of isometries $G$ of $M$ acts cocompactly on $M$, the strong
Sobolev inequality in $L^p$ is satisfied if and only if $G$ is
either non-amenable or non-unimodular \cite{Pit}.

Assume now that $M=G$ is a connected, unimodular Lie group, endowed
with a left invariant Riemannian metric. Denote by $\rho_{G,p}$ the
right regular representation on $D_p(G)$. Let $g\in G$ and
$\gamma:\;[0,d(1,g)]\to G$ be a geodesic between $1$ and $g$. For
any $f\in D_p(G)$ and $x\in G$, we have
$$(f-\rho_{G,p}(g)f)(x)=f(x)-f(xg)=\int_0^{d(1,g)}\nabla f(\gamma_x(t))\cdot\gamma_x'(t)dt,$$
where $\gamma_x(t)=x\gamma(t)$. Using H\"older's inequality, we
deduce that
$$\|f-\rho_{G,p}(g)f\|_p\leq d(1,g)\|\nabla f\|_p.$$
Therefore, there is a well defined map from $D_p(G)$ to
$Z^1(G,\rho_{G,p})$
$$J:\;f\mapsto (b_f:\; g\mapsto f-\rho_{G,p}(g)f).$$
The map $J$ induces an injective map from $\D_p(G)$ to
$Z^1(G,\rho_{G,p})$. Moreover, $b_f$ is a coboundary if and only if
$f$ is in $L^p(G)+\{\textnormal{constants}\}$, i.e. if and only if
the class of $f$ is zero in $H_p^1(G)$. Hence, $J$ induces an
injective\footnote{Actually, $J$ induces an isomorphism of
topological vector spaces \cite[Proposition~6.3]{Te} but this is
more delicate and not needed here.} linear map from $H_p^1(G)$ to
$H^{1}(G,\rho_{G,p})$.

Let $G$ be a simple  Lie group of rank 1 equipped with a
left-invariant Riemannian metric. Up to taking the quotient by a
normal compact subgroup, $G$ is $\textnormal{PO}(n,1)$,
$\textnormal{PU}(n,1)$, $\textnormal{PSp}(n,1)$ or
$\textnormal{F}_{4(-20)}$. Let $\partial G$ be the sphere at
infinity of $G$, and let $p_0$ be its Hausdorff dimension, so that
$$p_0\;=\;
\left\{
\begin{array}{ccc}
n-1 & if & G=PO(n,1)\\
2n & if & G=PU(n,1), n\geq 2\\
4n+2 & if & G=PSp(n,1)\\
22 & if & G=F_{4(-20)}
\end{array}
\right.  $$ By a result of P. Pansu \cite{Pan,Pansu},
$H^1_{(p)}(G)\neq 0$ if and only if $p>\max\{1,p_0\}$. From the
above discussion, we deduce that ${H}^1(G,\rho_{G,p})\neq 0$ for
those groups as soon as $p>\max\{1,p_0\}$. Together with the fact
that connected Lie groups have Property $(BP_0^{L^p})$ for
$1<p<\infty$, this yields the following result.

\begin{thm}\label{rank1}
Let $G$ be a simple Lie group of rank $1$ over $k$. Let $p_0$ be the
Hausdorff dimension of the visual boundary of $G$. Then, for every
$p>\max\{1,p_0\}$, there exists a proper affine action of $G$ on
$L^p(G)$ with linear\footnote{Since $G$ is unimodular, the
representations $\lambda_{G,p}$ and $\rho_{G,p}$ are isomorphic.}
part $\lambda_{G,p}$. \hfill$\qed$
\end{thm}

\section{Actions on Hilbert spaces}

\subsection{Property $({\BP}_0)$}

Recall that a locally compact group $G$ has Property (FH) if every
affine isometric action of $G$ on a Hilbert space has a fixed point.
For $G$ $\sigma$-compact, this is known to be equivalent to the
celebrated Kazhdan's Property~(T) (see \cite[Chapter~2]{BHV}).

When $V$ is a Hilbert space (sufficiently large in comparison to
$G$), we write $({\BP}_0)$ and {\FHO} for $(\BP_0^V)$ and
$(\FV_0^V)$.

There is a simple characterization of groups with Property {\FHO}
among groups with Property $({\BP}_0)$.

\begin{prop}
Let $G$ be a locally compact group.
\begin{enumerate}
\item [1)] Suppose that $G$ has Property
$({\BP}_0)$. Then either $G$ is a-T-menable or has Property {\FHO}.

\item [2)] If $G$ is both a-T-menable and has Property {\FHO}, then
it is compact.
\end{enumerate}\label{prop:FH0_haag}
\end{prop}
\begin{proof} The first statement is clear. Suppose that $G$ is
a-T-menable and is not compact. Then $G$ is $\sigma$-compact (take a
proper action $\alpha$ and consider $K_n=\{g\in
G\,:\;\|\alpha(g)(0)\|\le n\}$). Since $G$ is a-T-menable, it is
Haagerup, i.e. it has a $C_0$-representation $\pi$ with almost
invariant vectors; since $G$ is not compact, $\pi$ has no invariant
vector. By Proposition 2.5.3 in \cite{BHV}, $\infty\pi$ has
nontrivial 1-cohomology, while it is $C_0$. Hence $G$ does not have
Property {\FHO}.\end{proof}
\medskip

Let us mention an application of Property $({\BP}_0)$ in ergodic
theory.

\begin{prop}\label{prop:mixing} Let $G$ be a locally compact group
with Property $({\BP}_0)$ and $Hom(G,\R)=0$. Let $G$ act (on the
right), in a measure-preserving way, on a probability space
$(X,{\cal B},\mu)$; assume that the action is mixing. Let $F:X\times
G\rightarrow\C$ be a measurable function such that
\begin{itemize}
\item $\int_X |F(x,g)|^2\,d\mu(x)<\infty$ for every $g\in G$;
\item $F(x,gh)=F(x,g)+F(x.g,h)$ for every $g,h\in G$, almost
everywhere in $x$.
\end{itemize}
Then the following alternative holds: either there exists $f\in
L^2(X,\mu)$ such that $F(x,g)=f(x.g)-f(x)$ (for every $g\in G$,
almost everywhere in $x$), or $$\lim_{g\rightarrow\infty}\int_X
|F(x,g)|^2\,d\mu(x)=\infty.$$
\end{prop}

\begin{proof} Set $L^2_0(X,\mu)=\{f\in L^2(X,\mu):\int_X
f(x)\,d\mu(x)=0\}$. As $\int_X F(x,gh)\,d\mu(x)=\int_X
F(x,g)\,d\mu(x)+\int_X F(x,h)\,d\mu(x)$ for every $g,h\in G$, we
have $F(\cdot,g)\in L^2_0(X,\mu)$ for every $g\in G$, in view of the
assumption $\textnormal{Hom}(G,\R)=0$. Let $\pi$ be the standard
representation of $G$ on $L^2_0(X,\mu)$. Then $F(\cdot,g)$ defines a
1-cocycle with respect to $\pi$. Since the $G$-action on $X$ is
mixing, the representation $\pi$ is $C_0$ (see Theorem 2.9 in
\cite{BekMay}), so that the conclusion follows immediately from
Property $({\BP}_0)$.
\end{proof}

\subsection{Discrete groups without $({\BP}_0)$}

Proposition \ref{groupswithPH0} provides a wealth of groups with
Property $({\BP}_0)$. We now provide examples of groups {\it
without} Property $({\BP}_0)$ (in particular the free group $F_{n}$
on $n$ generators, $n\geq 2$).

\begin{prop}\label{amalgams}
Let $H$ be an infinite group, $K$ a non-trivial group, and $F$ a
common finite subgroup of $H$ and $K$, which is distinct from $K$.
Let $G=H\ast_{F}K$ be the amalgamated product. Then there exists a
1-cocycle with respect to the regular representation $\lambda_G$
which is neither bounded nor proper. In particular, $G$ does not
have Property $({\BP}_0)$.\label{prop:amalgam_non_PH0}
\end{prop}
\begin{proof} Let $w$ be a $\ell^2$ function on $G$ which is left
$F$-invariant, but not left $K$-invariant (in particular $w\neq 0$).
Define $\alpha(k)=\lambda_G(k)$ for $k\in K$, and
$\alpha(h)=t_w\circ \lambda_{G}(h)\circ t_{-w}$ for $h\in H$, where
$t_w$ denotes the translation by $w$ on $\ell^2(G)$. Then $\alpha$
is well-defined on $H\ast_F K$ (by the $F$-invariance assumption on
$w$). The fixed point set of $K$ is the set of all left
$K$-invariant functions. The set of fixed points of $H$ is reduced
to $\{w\}$ (since $H$ is infinite). Accordingly, the action has no
fixed point. On the other hand, since $H$ is infinite and has a
fixed point, the action is not proper.\end{proof}
\medskip

To obtain other examples of groups without \PHO, we first establish
a connection with a classical conjecture on discrete groups. For a
group $\Gamma$, we denote by $\C\Gamma$ the group algebra over $\C$,
and by denote again by $\lambda_{\Gamma}$ the left regular
representation of $\C\Gamma$ on $\ell^2(\Gamma)$:
$$\lambda_{\Gamma}(f)\xi = f\ast\xi$$
($f\in\C\Gamma,\,\xi\in\ell^2(\Gamma)$). Here is the {\it analytical
zero-divisor conjecture}:

\begin{conj}\label{zerodiv} If $\Gamma$ is a torsion-free group,
then $\lambda_{\Gamma}(f)$ is injective, for every
non-zero $f\in\C\Gamma$.
\end{conj}

The main result on Conjecture \ref{zerodiv} is due to P. Linnell
\cite{Linn}: it holds for groups which are extensions of a
right-orderable group by an elementary amenable group; in
particular, we will use the fact that it holds for non-abelian free
groups.

\begin{lem}\label{aagard} Let $\Gamma$ be a group satisfying Conjecture
\ref{zerodiv}. Let $f_1,f_2 \in\C\Gamma$ be non-zero
elements. There exists non-zero functions
$\xi_1,\xi_2\in\ell^2(\Gamma)$ such that $\lambda_{\Gamma}(f_1)\xi_1
+ \lambda_{\Gamma}(f_2)\xi_2 \,=\,0$.
\end{lem}

\begin{proof} We start with a

{\bf Claim:} If $f\in\C\Gamma$ is a non-zero element, then
$\lambda_{\Gamma}(f)$ has dense image. To see that, observe that the
orthogonal of the image of $\lambda_{\Gamma}(f)$ is the kernel of
$\lambda_{\Gamma}(f^*)$, which is $\{0\}$ as $\Gamma$ satisfies
Conjecture \ref{zerodiv}.
\medskip

Let $L(\Gamma)$ be the von Neumann algebra of $\Gamma$, i.e. the
bi-commutant of $\lambda_{\Gamma}(\C\Gamma)$ in ${\cal
B}(\ell^2(\Gamma))$. A (non-necessarily closed) subspace of
$\ell^2(\Gamma)$ is {\it affiliated with $L(\Gamma)$} if it is
invariant under the commutant $\lambda_{\Gamma}(\C\Gamma)'$ of
$\lambda_{\Gamma}(\C\Gamma)$. E.g., if $f\in\C\Gamma$, the image of
$\lambda_{\Gamma}(f)$ is an affiliated subspace. A result of L.
Aagaard \cite{Aag} states that the intersection of two dense,
affiliated subspaces is still dense. We apply this with the images
of $\lambda_{\Gamma}(f_1)$ and of $\lambda_{\Gamma}(-f_2)$, so that
there exist non-zero $\xi_1,\xi_2$ such that
$\lambda_{\Gamma}(f_1)\xi_1 = \lambda_{\Gamma}(-f_2)\xi_2$.
\end{proof}

\begin{prop}\label{coconfreegrps}
Fix $k\geq 2$. Let $w$ be a non-trivial reduced word in the free
group $F_k$. There exists an unbounded 1-cocycle $b_w\in
Z^1(F_k,\lambda_{F_k})$, such that $b_w(w)=0$. In particular, $b_w$
is not proper.
\end{prop}

\begin{proof} We start with $k=2$.  Write $F_2$ as the free group
on 2 generators $s,t$. Write $w$ as a reduced word in $s^{\pm 1},
t^{\pm 1}$:
$$w\,=\,x_1^{\epsilon_1}x_2^{\epsilon_2}...x_n^{\epsilon_n}$$
($x_j\in\{s,t\};\epsilon_j \in\{-1,1\}$). If either $s$ or $t$ does
not appear in $w$, then the existence of the desired cocycle follows
from the proof of Proposition \ref{amalgams} (with $H=K=\Z$). So may
assume that both $s$ and $t$ appear in $w$. Set $\delta_j =
\frac{\epsilon_j -1}{2}$ and define two elements
$f_{w,s},\,f_{w,t}\in \C F_2$ by:
$$f_{w,s}\,=\,\sum_{j: x_j=s}\epsilon_j x_1^{\epsilon_1}...x_{j-1}^{\epsilon_{j-1}}x_j^{\delta_j};$$
$$f_{w,t}\,=\,\sum_{j: x_j=t}\epsilon_j x_1^{\epsilon_1}...x_{j-1}^{\epsilon_{j-1}}x_j^{\delta_j}.$$
Note that $f_{w,s}$ and $f_{w,t}$ are non-zero, as $s$ and $t$
appear in $w$. Since $F_2$ satisfies Conjecture \ref{zerodiv} (by
Linnell's result already quoted \cite{Linn}), we may appeal to Lemma
\ref{aagard} and find non-zero functions $\xi_s,\xi_t\in\ell^2(F_2)$
such that $\lambda_{F_2}(f_{w,s})\xi_s + \lambda_{F_2}(f_{w,t})\xi_t
\,=\,0$.

Set then $b_w(s)=\xi_s,\,b_w(t)=\xi_t$ and, using freeness of $F_2$,
extend uniquely to a 1-cocycle $b_w\in Z^1(F_2,\lambda_{F_2})$.
Using the relations $$b(g_1g_2...g_m)=\sum_{j=1}^m
\lambda_{F_2}(g_1...g_{j-1})b(g_j)$$ and
$b(g^{-1})=-\lambda_{F_2}(g^{-1})b(g)$ (for a cocycle $b$ and
$g_1,...,g_m,g \in F_2$), one checks that
$$b_w(w)\,=\,\lambda_{F_2}(f_{w,s})b_w(s)+\lambda_{F_2}(f_{w,t})b_w(t)\,=\,0.$$
It remains to show that $b_w$ is unbounded, i.e. that the
corresponding affine action $\alpha_w$ has no fixed point. Let
$H=<w>$ be the cyclic subgroup generated by $w$. As the linear
action is $C_0$, the only fixed point of $\alpha_w|_H$ is $0$. But
$0$ is clearly not fixed under $\alpha_w(s)$ or $\alpha_w(t)$, which
completes the proof in case $k=2$.
\medskip

Suppose now $k\geq 2$. View $F_k$ as a subgroup of index $k-1$ in
$F_2$. The restriction of $\lambda_{F_2}$ to $F_k$ is the direct sum
of $k-1$ copies of $\lambda_{F_k}$. Project the cocycle $b_w$ given
by the first part of the proof, to each of these $k-1$ summands.
This way, get $k-1$ cocycles in $Z^1(F_k,\lambda_{F_k})$, each of
them vanishing on $w$. At least one of them is unbounded, because
$b|_{F_k}$ is unbounded.
\end{proof}

\begin{cor}\label{amalgoverZ} For $k\geq 2$, let $\Gamma =F_k\ast_{\Z}G$ be
an amalgamated product over $\Z$ an
infinite cyclic subgroup. Then $\Gamma$ does not have property \PHO.
\end{cor}

\begin{proof} Let $w\in F_k$ and $g\in G$ correspond to the positive
generators of the copies of $\Z$ that are
amalgamated. Choosing representatives for the left cosets of $F_k$
in $\Gamma$, identify $\lambda_{\Gamma}|_{F_k}$ with
$\infty\lambda_{F_k}=:\lambda_{F_k}\oplus\lambda_{F_k}\oplus...$.
Let $b_w\in Z^1(F_k,\lambda_{F_k})$, as in Proposition
\ref{coconfreegrps}. Define an affine action $\alpha$ of $F_k$, with
linear part $\lambda_{\Gamma}|_{F_k}$, by:
$$\alpha(x)(v_1,v_2,v_3,...)\,=\,(\lambda_{F_k}(x)v_1 +b_w(x),\lambda_{F_k}(x)v_2,\lambda_{F_k}(x)v_3,...)$$
$(x\in F_k)$. On the other hand, view $\lambda_{\Gamma}|_G$ as an
affine action of $G$. Since
$$\alpha(w)=\lambda_{\Gamma}(w)=\lambda_{\Gamma}(g),$$
these two affine actions can be "glued together", i.e. extend to an
affine action $\tilde{\alpha}$ of $\Gamma$, with linear part
$\lambda_{\Gamma}$. By the very construction, $\tilde{\alpha}$ has
unbounded orbits and is not proper.
\end{proof}

\begin{cor}\label{surfgrps} The surface groups $\Gamma_g\;(g\geq 2)$ do not have Property \PHO.
\end{cor}

\begin{proof} The presentation
$$\Gamma_g\,=\,\langle a_1,...,a_g,b_1,...,b_g|[a_1,b_1]^{-1}=\prod_{j=2}^g [a_j,b_j]\rangle$$
exhibits $\Gamma_g$ as an amalgamated product $F_2
\ast_{\Z}F_{2g-2}$ so Corollary \ref{amalgoverZ} applies.
\end{proof}

Here is an intriguing question, in view of the fact that $PSL_2(\Z)$
contains a free group of finite index:

\begin{que} Does $PSL_2(\Z)\simeq C_2\ast C_3$ have Property
$({\BP}_0)$?\end{que}

\subsection{Application to the regular representation}

Let us recall that Guichardet \cite[Th\'eor\`eme~1]{Gui72} proved
that, if $\pi$ is a representation without non-zero fixed vector of
a locally compact, $\sigma$-compact group, the space $B^1(G,\pi)$ is
closed in $Z^1(G,\pi)$ (endowed with the topology of uniform
convergence on compact subsets) if and only if $\pi$ does not almost
have invariant vectors. In particular $H^1(G,\pi)\neq 0$ if $\pi$
almost has invariant vectors. This rests on a clever use of the open
mapping theorem for Fr\'echet spaces. Using this, we can reprove the
following result, first proved in \cite{AW} (see also
\cite{BekCheVal}).

\begin{prop}
Let $G$ be a $\sigma$-compact, locally compact group. If $G$ is
Haagerup, then it is a-T-menable.\label{prop:Haag_aTmenable}
\end{prop}
\begin{proof} Set $H=G\times\mathbf{Z}$; then $H$ is $\sigma$-compact,
locally compact, is Haagerup, and has noncompact center. Hence, by
Proposition \ref{noncompactcentre}, it has Property $({\BP}_0)$.
Take a $C_{0}$-representation $\pi$ of $H$, almost having invariant
vectors. By Guichardet's result recalled above, there exists an
affine action $\alpha$ of $H$, with linear part $\pi$, and without
fixed point. By Property $({\BP}_0)$, the action $\alpha$ is proper.
So the restriction $\alpha|_{G}$ is proper too.\end{proof}

If $G$ is $\sigma$-compact and amenable, the representation $\pi$ in
the above proof can be taken as the left regular representation of
$G\times\mathbf{Z}$ on $L^2(G\times\mathbf{Z})$. (By way of
contrast, if $\Gamma$ is a discrete, non-amenable group, then
$H^1(\Gamma\times\mathbf{Z},\lambda_{\Gamma\times\mathbf{Z}}) =0$ by
Corollary~10 in \cite{BekVal}).

Concerning affine actions on $L^2(G)$, we have the following

\begin{conj}
For an amenable group $G$, every affine action with linear part
$\lambda_{G}$ is either bounded or proper.
\end{conj}

Evidence for this conjecture comes from the fact that Proposition
\ref{PH0heredit}, Corollary \ref{solvable} and Proposition
\ref{groupswithPH0} establish it in numerous cases: amenable groups
with infinite center, solvable groups, amenable Lie groups,
etc\ldots More evidence comes from a result proved in \cite{MarVal}:
if $\Gamma$ is a countable amenable group, and $A$ is any infinite
subgroup, then the restriction map
$H^1(\Gamma,\lambda_{\Gamma})\rightarrow
H^1(A,\lambda_{\Gamma}|_{A})$ is injective. If true, our conjecture
would provide a conceptual explanation of this fact.

Being more ambitious, one may even ask

\begin{que} Does every amenable group have Property $({\BP}_0)$?
\end{que}

We now turn to the study of some groups $G$ for which
$H^1(G,\lambda_G)\neq 0$.

\begin{lem}\label{lem:vanish_1coho_product}
Let $G$ be a locally compact, second countable group. Suppose that,
for some $k\ge 2$, the group $G$ has closed normal subgroups
$N_1,\dots,N_k$ such that $[N_i,N_j]=1$ whenever $i\neq j$ and
$G=N_1\dots N_k$. Let $\pi$ be a unitary representation such that
$\overline{H^1}(G,\pi)\neq 0$. Then at least one of the $N_i$ has an
invariant vector by $\pi$.
\end{lem}
\begin{proof}
There is an obvious map $p$ of $N=\prod_{i=1}^kN_i$ onto $G$. Then
$\overline{H^1}(N,\pi\circ p)\neq 0$. This uses the standard fact
that every compact subset of $G$ is the image of a compact subset of
$N$ (note that we use here $\sigma$-compactness).

Suppose that for some $i$, the group $N_i$ has no invariant vector
by $\pi\circ p$. Write $N=N_i\times\prod_{j\neq i}N_j$; by
\cite[Proposition 3.2]{ShaInv} (which uses second countability),
$\prod_{j\neq i}N_j$ has an invariant vector by $\pi\circ p$, so
that for every $j\neq i$, $N_j$ has an invariant vector by $\pi\circ
p$.
\end{proof}

\begin{prop}\label{prop:H1bar_C0}
Let $G$ be a connected Lie group or $G=\mathbf{G}(\mathbf{K})$, the
group of $\mathbf{K}$-points of a linear algebraic group
$\mathbf{G}$ over a local field $\mathbf{K}$ of characteristic zero.
Suppose that $G$ has a $C_0$-representation $\pi$ such that
$\overline{H^1}(G,\pi)\neq 0$. Then either $G$ is amenable, or has a
compact subgroup $K$ such that $G/K$ is a simple Lie group (resp. a
simple linear algebraic group) with trivial centre.
\end{prop}
\begin{proof}By Property $({\BP}_0)$, $G$ has the Haagerup Property.
If $G$ is a connected Lie group, by
\cite[Chap.~4]{CCJJV}, $G=RS_1\dots S_k$ where $R,S_1,\dots,S_k$
centralize each other, $R$ is a connected amenable Lie group, and
each $S_i$ is a simple, noncompact, connected Lie group with the
Haagerup Property (with possibly infinite centre). In the case of an
algebraic group, the same conclusion holds \cite{CorJLT}, except
that the $S_i$'s are simple linear algebraic groups.

If $G$ is not amenable, then $k\ge 1$, and in this case by Lemma
\ref{lem:vanish_1coho_product} it follows that $k=1$ and $R$ is
compact. By
    \cite[Corollary~3.6]{ShaInv}, the centre $Z(G)$ has an
    invariant vector by $\pi$ and thus is compact since $\pi$ is $C_0$;
    since in our situation $Z(S_1)\subset Z(G)$, we see that
    $S_1$ has finite centre, so that $K=RZ(S_1)$ is compact and
    $G/K$ is a simple group with trivial centre.
\end{proof}

\begin{prop}\label{H1reg}
Let $G$ be a connected Lie group or $G=\mathbf{G}(\mathbf{K})$ where
$\mathbf{K}$ is a local field of characteristic zero. Assume $G$
non-compact. Then the following are equivalent\begin{itemize}
    \item[(i)] $H^1(G,\lambda_G)\neq 0$.
    \item[(ii)] Either $G$ is amenable, or there exists a compact normal subgroup $K\subset G$ such
    that $G/K$ is isomorphic to $\textnormal{PSL}_2(\R)$ (case of Lie groups),
    or a simple algebraic group of rank one (case of an algebraic group over a $p$-adic field).
    \end{itemize}
\end{prop}

\begin{proof}
Suppose (i). If $G$ is not amenable, then, by the result of
Guichardet already mentioned \cite[Th\'eor\`eme~1]{Gui72}, one has
$\overline{H^1}(G,\lambda_G)=H^1(G,\lambda_G)\neq 0$. By Proposition
\ref{prop:H1bar_C0}, $G$ has a compact normal subgroup $K$ such that
$S=G/K$ is simple with trivial centre. Moreover, $G$ does not have
Property~(T), hence has rank one \cite{DK}. This settles the
non-Archimedean case. If $G$ is a Lie group, then by
\cite[Theorem~6.4]{Mar}, $\lambda_G$ contains an irreducible
subrepresentation $\sigma$ factoring through $S$, such that
$\overline{H^1}(G,\sigma)=\overline{H^1}(S,\sigma)\neq 0$. Then
$\sigma\le\lambda_S$, as $S$ is cocompact in $G$, so that
$H^1(S,\lambda_S)\neq 0$. By a result of Guichardet (Proposition 8.5
in Chapter III of \cite{Guich}), this implies that
$S\simeq\textnormal{PSL}_2(\R)$.

Conversely suppose (ii). If $G$ is amenable, then $H^1(G,\lambda_G)$
is not Hausdorff, hence is nonzero. Otherwise, suppose $G$
non-amenable, and consider $K$ as in (ii). By Proposition
\ref{trees} and Theorem \ref{rank1} (noticing that $p_0=1$ for
$PSL_2(\R)$), we have $H^1(G/K,\lambda_{G/K})\neq 0$. Then
$H^1(G,\lambda_G)\neq 0$ by the same elementary argument as used in
the proof of Proposition \ref{trees}. \end{proof}

\begin{cor}\label{bettiL2} Let $G$ be a connected Lie group or $G=\mathbf{G}(\mathbf{K})$
where $\mathbf{K}$ is a local field of characteristic zero; let
$\Gamma$ be a uniform lattice in $G$. If the first $L^2$-Betti
number $\beta^1_{(2)}(\Gamma)$ is non-zero, then $\Gamma$ is
commensurable either to a non-abelian free group or to a surface
group (more precisely: $\Gamma$ has a finite index subgroup
$\Gamma_0$ with a finite normal subgroup $N$ such that $\Gamma_0/N$
is either a non-abelian free group or a surface group).
\end{cor}

\begin{proof} From $\beta^1_{(2)}(\Gamma)>0$, it follows that $\Gamma$
(and also $G$) is non-amenable: see Theorem 0.2 in \cite{CheGro}. On the other
hand, it was proved in \cite{BekVal} that, for $\Gamma$ a finitely
generated non-amenable group: $\beta^1_{(2)}(\Gamma)>0$ if and only
if $H^1(\Gamma,\lambda_{\Gamma})\neq 0$. Since $\Gamma$ is uniform
in $G$, we have by Shapiro's lemma (Proposition 4.6 in Chapter III
of \cite{Guich}):
$$0\neq H^1(\Gamma,\lambda_{\Gamma})=
H^1(G,\textnormal{Ind}^G_{\Gamma}\lambda_{\Gamma})\simeq H^1(G,\lambda_G).$$
By Proposition \ref{H1reg}, the group $G$ admits a compact normal
subgroup $K$ such that $G/K$ is isomorphic either to
$\textnormal{PSL}_2(\R)$ or to a simple algebraic group of rank 1
over a $p$-adic group. Let $p:G\rightarrow G/K$ be the quotient map.
Then $p(\Gamma)$ is a uniform lattice in $G/K$. By Selberg's lemma,
find a finite-index torsion-free subgroup $\tilde{\Gamma}_0$ of
$p(\Gamma)$: then $\tilde{\Gamma}_0$ is either a surface group (case
of $\textnormal{PSL}_2(\R)$) or a non-abelian free group
(non-Archimedean case). Set $\Gamma_0 =p^{-1}(\tilde{\Gamma}_0)$, a
finite-index subgroup of $\Gamma$. Conclude by observing that the
kernel $\Ker (p|_{\Gamma_0})$ is contained in $\Gamma\cap K$, so is
finite.
\end{proof}

The preceding result overlaps a result of B. Eckmann (Theorem 4.1 in
\cite{Eck}), who classified lattices $\Gamma$ (not necessarily
uniform) with $\beta^1_{(2)}(\Gamma)>0$ in a connected Lie group.

\subsection{Some non-$\sigma$-compact groups}

Here is a curiosity. Start with the observation from the proof of
Proposition \ref{prop:FH0_haag} that a locally compact,
non-$\sigma$-compact group cannot be a-T-menable. Accordingly, if it
has Property $({\BP}_0)$, then it also has Property {\FHO}.

The above observation shows that the $\sigma$-compactness assumption
is necessary in Guichardet's result mentioned above. It also
provides, in the non-Fr\'echet case, some explicit counterexamples
to the statement of the open mapping theorem.

\begin{prop} Let $G$ be a non-$\sigma$-compact locally compact amenable
group with Property $({\BP}_0)$. Endow $Z^1(G,\lambda_{G})$ with the
topology of uniform convergence on compact subsets. Then the map
$$\partial:
\left\{
\begin{array}{ccc}
L^2(G) & \rightarrow & Z^1(G,\lambda_{G}) \\
\xi & \mapsto & (g\mapsto \lambda_{G}(g)\xi-\xi)
\end{array}\right.$$ is a continuous linear
bijective homomorphism, whose inverse is not continuous.
\end{prop}

\begin{proof} The map $\partial$ is linear, injective (as $G$ is not
compact) and surjective (since $G$ has property ({\FHO})). It is
clearly continuous. By amenability of $G$, for every $\eps>0$ and
every compact subset $K\subset G$, there exists a unit vector
$\xi\in L^2(G)$ such that
$$\max_{g\in K}\|(\partial\xi)(g)\| < \eps.$$
This clearly shows that $\partial^{-1}$ is not
continuous.\end{proof}

\begin{exe} Examples of non-$\sigma$-compact amenable groups with Property
$({\BP}_0)$ include
\begin{itemize}
    \item  Uncountable solvable groups (by Corollary \ref{solvable})

    \item  Discrete groups of the form $G=F^{I}$, where $F$ is a
non-trivial
    finite group and $I$ is any infinite set. Indeed $G$ is amenable, as
    it is locally finite, and since $G$ is isomorphic to $G\times G$
    (as $I$ is infinite), $G$ contains an infinite normal subgroup with
    infinite centralizer, so Proposition \ref{PH0heredit} applies.
\end{itemize}
\end{exe}

\section{Actions of $\mathbf{Z}$ and $\mathbf{R}$}

\subsection{Actions of $\mathbf{Z}$}

We have shown that every action of $\Z$ on a Hilbert space with
$C_0$ linear part is either bounded or proper.

An example of Edelstein \cite{Edel} shows that the $C_0$ assumption
cannot be dropped. Let us briefly recall his example. On $\C$,
consider the rotation $r_n$ with centre $1$ and angle $2\pi/n!$.
Consider the abstract product $\C^\N$, and, for
$(z_n)_{n\in\N}\in\C^\N$, set
$r((z_n)_{n\in\N})=(r_n(z_n))_{n\in\N}$. This self-map is bijective
and has the constant sequence 1 as unique fixed point. Moreover, it
can be shown that $r(\ell^2(\N))=\ell^2(\N)$. Thus $r$ induces an
affine isometry of $\ell^2(\N)$, which has no fixed point since the
constant 1 is not in $\ell^2(\N)$. However, the action is not
proper; actually 0 is a recurrent point: an easy computation gives
$\|r^{n!}(0)\|^2\le\sum_{k>n}(2\pi n!/k!)^2$, and this sum clearly
tends to zero. Notice that $r$ almost has fixed points: with $v_m$
the characteristic function of $\{1,...,m\}$, one has
$\lim_{m\rightarrow\infty}\|r(v_m)-v_m\|=0$.

Observe that this isometry has diagonalizable linear part. Let us
now provide another counter-example with further assumptions on the
linear part.

\begin{defn}
A unitary or orthogonal representation of a group is \textit{weakly
$C_0$} if it has no nonzero finite dimensional
subrepresentation\footnote{$C_0$ (resp. weakly $C_0$)
representations are often called mixing (resp. weakly mixing).}.
\end{defn}

\begin{prop}
There exists an affine isometric action of $\Z$ on a complex Hilbert
space, which is neither bounded nor proper, and has weakly $C_0$
linear part.\label{prop:affine_Z_non_wc0}
\end{prop}

\begin{proof} Write $\sigma$ for the affine action of $\mathbf{Z}$, and $\pi$ for its linear
part. Let $\mu$ be a probability measure on $[0,1]$ and write
$H=L^2([0,1],\mu)$. Let $\pi(1)$ be the multiplication by the
function $e(x)=\exp(2i\pi x)$. Write $\sigma(1)=\tau_{1}\circ\pi(1)$
where $\tau_{1}$ is the translation by the constant function $1$.
Note that $\pi$ is weakly $C^0$ if and only if the spectrum of
$\pi(1)$ has no atom, i.e. $\mu$ is nonatomic.

Let $b$ be the corresponding cocycle and write $c(n)=\| b(n)\|^2$.
An immediate computation shows that
$$c(n)=\int\phi_n(x)d\mu(x)$$
where $\phi_n(x)=|\sin(\pi nx)/\sin (\pi x)|^2.$

Let $N_n$ a increasing sequence of integers and let $\eps_n$ be a
decreasing sequence in $]0,1[$, such that $\eps_n \rightarrow 0$.
Moreover, let us assume that $N_n/N_{n+1}=o(\eps_n)$.

For all positive integer $n$, write $$I_n(k)=
\left[\frac{k-\eps_n}{N_n}, \frac{k+\eps_n}{N_n}\right]\cap [0,1]$$
and
$$K_n= K_{n-1}\cap \bigcup_{k=0}^{N_n} I_n(k).$$ Finally, write
$$K=\bigcap_n
K_n.$$ One can check easily that $K$ is homeomorphic to a Cantor
space.

Let $\mu$  be a probability measure on $[0,1]$ such that
\begin{itemize}
\item its support is contained in $K\cap [0,1/2]$; \item There
exists a subsequence $\eps_{k_n}$ such that
\begin{equation}\label{p}
\mu\left([0,\sqrt{\eps_{k_n}}]\right)=\sqrt{\eps_{k_n}}.
\end{equation}
\end{itemize}
We choose the sequence $k_n$ such that each interval
$I_n=[\sqrt{\eps_{k_{n+1}}}, \sqrt{\eps_{k_n}}]$ intersects $K$
nontrivially. Take for $\mu_n$ any nonatomic measure supported by
$K\cap I_n$ such that $\mu_n(I_n)=\sqrt{\eps_{k_{n}}}-
\sqrt{\eps_{k_{n+1}}}$ and define $\mu=\sum_n \mu_n:$ clearly, $\mu$
is nonatomic.

{\bf Claim 1.} The action $\sigma$ has no fixed point (so it has
unbounded orbits).

If $\sigma$ has a fixed point $f$, then $f(x)=(1-\exp(2i\pi
x))^{-1}$ $\mu$-a.e. Let us show that $f$ does not belong to
$L^2([0,1])$. Indeed, note that $|f|^2=\left(1/\sin(\pi
x)\right)^2$. For all $x\in [0,\sqrt{\eps_{k_n}}]$, we have
$$\sin(\pi x)^2\leq \pi^2 x^2\leq \pi^2\eps_{k_n}$$
and by (\ref{p})
$$\mu([0,\sqrt{\eps_{k_n}}])= \sqrt{\eps_{k_n}}.$$
It follows that
$$\int|f|^2d\mu\geq \frac{1}{\pi^2\sqrt{\eps_{k_n}}}$$
which proves claim 1.

{\bf Claim 2.} If moreover $\eps_{k_n}=o(N_{k_n}^{-4})$ (for
instance, $N_n=2^{n!}$ and $\eps_n=(N_n)^{-5}$), then $c(N_{k_n})$
tends to $0$, so that the action is not proper.

Indeed, let us show that $c(N_{k_n})=o(1)$.

First, note that for all $x\in K$, the fractional part of
$N_{k_n}.x$ is less than $\eps_{k_n}$. Thus, for every $x\geq
\sqrt{\eps_{k_n}}$ and every $x\in K$, it comes
$$\phi_{N_{k_n}}(x)\leq \left(\frac{\sin(2\pi \eps_{k_n})}{\sin
(2\pi\sqrt{\eps_{k_n}} )}\right)^2\leq \pi^2\eps_{k_n}/4.$$ On the
other hand, we have
$$\frac{\sin (2\pi N_{k_n}x)}{\sin (2\pi x)} \leq N_{k_n}$$ and by
(\ref{p})
$$\mu([0, \sqrt{\eps_{k_n}}])=\sqrt{\eps_{k_n}}.$$ It follows that
$$c(N_{k_n})\leq \sqrt{\eps_{k_n}}.N_{k_n}^2+\pi^2\eps_{k_n}/4.$$
So we get $c(N_{k_n})=o(1).$
\end{proof}

\subsection{Actions of $\mathbf{R}$}

Let us now show that the ``pathological" actions of $\Z$ described
above can be extended to $\R$.

Recall that a group $G$ is said to be {\it exponential} if, for
every $g\in G$, there is a one-parameter subgroup through $g$ (i.e.
a continuous homomorphism $\beta:\R\rightarrow G$ such that
$\beta(1)=g$). Clearly, an exponential group has to be
arc-connected.

Endow the group of affine isometries of a complex Hilbert space,
$\HH\rtimes\mathcal{U}(\mathcal{H})$, with the product topology, for
the natural topology on $\HH$ and the norm operator topology on the
unitary group $\mathcal{U}(\mathcal{H})$.

\begin{prop}\label{Isomexpon} The group of affine isometries of a complex Hilbert
space $\HH$, is exponential.
\end{prop}

\begin{proof} Let $\alpha(v)=Uv + b$ be an affine isometry of $\HH$.

By the spectral theorem, we find a projection-valued measure
$\mathbf{P}$ on $[-\pi,\pi[$ such that
$U=\int_{-\pi}^{\pi}e^{ix}\,d\mathbf{P}(x)$, in the sense that, for
every $\xi,\eta\in\HH$
$$\langle
U\xi|\eta\rangle\,=\,\int_{-\pi}^{\pi}e^{ix}\,d\mu_{\xi,\eta}(x)$$
where $\mu_{\xi,\eta}(A)=\langle\mathbf{P}(A)\xi|\eta\rangle$ for
any Borel subset $A\subset [-\pi,\pi[$. Consider the one-parameter
group of unitary operators
$$\upsilon(s)=\int_{[-\pi,\pi[}e^{isx}\,d\mathbf{P}(x).$$

For every $\xi\in\HH$ and $t\in\R$, define
$$b_\xi(t)=\int_{0}^{t}\upsilon(s)\xi\,ds.$$
It is straightforward that $b_\xi\in Z^1(\R,\upsilon)$. Let us
consider the operator $A=\int_0^1 \upsilon(s)\,ds$. Then
$b_\xi(1)=A\xi$ for every $\xi\in\HH$. Thus, to show that there
exists $\xi\in\HH$ such that $b=b_{\xi}(1)$, it suffices to
establish that $A$ is invertible.

By Fubini's Theorem
$$\int_0^1 \upsilon(s)\,ds\;=\;\int_0^1 (\int_{[-\pi,\pi[}e^{isx}\,d\mathbf{P}(x))\,ds\;=\;
\int_{[-\pi,\pi[}\frac{e^{ix}-1}{ix}\,d\mathbf{P}(x).$$ Since the
function $x\mapsto \frac{ix}{e^{ix}-1}$ is bounded on $[-\pi,\pi[$,
we obtain that
$$\int_{[-\pi,\pi[}\frac{ix}{e^{ix}-1}\,d\mathbf{P}(x)$$
is a bounded operator on $\HH$, and is the inverse of $A$; so we may
take $\xi=A^{-1}(b)$.
\end{proof}

In view of Proposition \ref{prop:affine_Z_non_wc0}, we obtain

\begin{cor}
There exists an affine isometric action of $\R$ on a complex Hilbert
space that is neither bounded nor proper. Moreover, it can be chosen
weakly~$C_0$.\qed
\end{cor}

\begin{rem} Proposition \ref{Isomexpon} is {\it false} for {\it real} Hilbert
spaces. This follows from the fact that the orthogonal group of a
real Hilbert space is {\it not} exponential. This is clear in finite
dimension (the group $O(n)$ is not connected), and was observed by
Putnam and Wintner \cite{PuWi} in infinite dimension (although the
orthogonal group is then connected \cite{PuWi1}). An example of an
orthogonal transformation which is not in the image of the
exponential map is a reflection
$$S\,=\,\textnormal{diag}(-1,1,1,1,\ldots);$$
this can be seen by noticing that $S$ is not a square in the
orthogonal group: indeed if $S=R^2$, since $R$ commutes with $S$ it
stabilizes the $-1$-eigenspace of $S$, which leads to a
contradiction.
\end{rem}

\bigskip
\footnotesize
\noindent Yves de Cornulier\\
IRMAR, Campus de Beaulieu, \\
35042 Rennes Cedex, France\\
E-mail: \url{decornul@clipper.ens.fr}\\
\medskip

\noindent Romain Tessera\\
Department of mathematics, Vanderbilt University,\\
Stevenson Center, Nashville, TN 37240 USA,
\\E-mail: \url{tessera@clipper.ens.fr}
\medskip

\noindent Alain Valette\\
Institut de Mathématiques - Université de Neuchâtel\\
11, Rue Emile Argand - BP 158, 2009 Neuch\^atel - Switzerland\\
E-mail: \url{alain.valette@unine.ch}

\end{document}